\documentclass[12pt]{article}
\usepackage{a4,amsmath,amsthm,amssymb,amsfonts,hyperref,enumerate}

\theoremstyle{plain}
\newtheorem{thm}{Theorem}      

\newtheorem{lem}{Lemma}  
\newtheorem{prop}{Proposition}

\theoremstyle{definition}
\newtheorem{rem}{Remark}  
\newtheorem{conj}{Conjecture}
\newtheorem{prob}{Problem}

\newcommand{\C}{{\mathbb C}}

\newcommand{\abs}[1]{{\left| {#1} \right|}}
\newcommand{\p}[1]{{\left( {#1} \right)}}
\newcommand{\Oh}[1]{{O \p{#1}}}

\newcommand{\norm}[1]{\left \Vert {#1} \right \Vert}

\renewcommand{\Re}{\operatorname{Re}}

\author{Johan Andersson\footnote{Email:johana@math.su.se}}

\begin{document}

\title{Tur{\'an}'s problem 10 revisited}

\maketitle

\begin{abstract}
  In this paper we prove that
  \begin{gather*}
    \inf_{\abs{z_k} \geq 1} \max_{\nu=1,\ldots,n^2} \abs{\sum_{k=1}^n
      z_k^\nu} = \sqrt n+\Oh{n^{0.2625+\epsilon}}. \qquad (\epsilon>0)
  \end{gather*}
  This improves on the bound
  \begin{gather*}
    \inf_{\abs{z_k} \geq 1} \max_{\nu=1,\ldots,n^2} \abs{\sum_{k=1}^n
      z_k^\nu} \leq \sqrt{6 n \log(1+n^2)}
  \end{gather*}
  of Erd{\H o}s and Renyi.  In the special case of $n+1$ being a prime
  we have previously obtained the much sharper result
  \begin{gather*}
    \sqrt n \leq \inf_{\abs{z_k} \geq 1} \max_{\nu=1,\ldots,n^2}
    \abs{\sum_{k=1}^n z_k^\nu} \leq \sqrt{n+1}.
  \end{gather*}
  The method of proof combines a general lower bound (of Andersson),
  explicit arithmetical constructions (of Montgomery, Fabrykowski or
  Andersson), moments (probabilistic methods) and estimates for the
  difference of consecutive primes (of
  Baker-Harman-Pintz). We also prove some (conditional and unconditional) related results. \end{abstract}

\section{Introduction}

In his book on the power sum method \cite{Turan} Tur{\'a}n proposed as
problem 10 (page 190) the study of the power sum quantity
\begin{gather*}
  \inf_{\abs{z_k} \geq 1}\max_{\nu=1,\ldots,n^2} \abs{\sum_{k=1}^n
    z_k^\nu}.
\end{gather*}
In a previous paper \cite{Andersson} we proved the strong inequality
\begin{gather} \label{ii} \sqrt n \leq \inf_{\abs{z_k} \geq 1}
  \max_{\nu=1,\ldots,n^2} \abs{\sum_{k=1}^n z_k^\nu} \leq \sqrt{n+1}
\end{gather}
whenever $n+1$ is prime. A natural question to ask
is (Alexei Venkov asked us this question when we visited
  Aarhus): What about general integers $n$? This problem makes sense
for all positive integers $n$ and there seems to be nothing special a
priori with the primes.

The problem has turned out to be much more difficult for a general
integer $n$.  Erd{\H o}s and Renyi \cite{Erdos} proved that
\begin{gather} \label{erdren} \inf_{\abs{z_k} \geq 1}
  \max_{\nu=1,\ldots,n^2} \abs{\sum_{k=1}^n z_k^\nu}=\Oh{ \sqrt{n \log
      n}},
\end{gather}
which follows from the following more general Lemma.
\begin{lem} {\rm (Erd{\H o}s and R{\'e}nyi)} \label{erd} There exists
  an n-tuple $(z_1,\ldots,z_n)$ of complex numbers such that $|z_k|=1$
  and $\max_{\nu=1,\dots,m} |\sum_{k=1}^n z_k^\nu| \leq \sqrt {6 n
    \log (m+1)}$.
\end{lem}
\noindent This remains the best bound so far for a general integer $m$. The proof of the Erd\H os-Renyi's lemma uses probabilistic methods and  is non constructive. Tijdeman (\cite{Tijdeman} and Tur\'an \cite{Turan} page 82) has given an explicit construction which implies the Erd{\H o}s-R{\'e}nyi lemma 
for $m=n^A$ (with a constant depending on $A>1$). 

 When we wrote our paper \cite{Andersson} we had no idea even how to prove
that
\begin{gather} \label{ajk} \inf_{\abs{z_k} \geq 1}
  \max_{\nu=1,\ldots,n^2} \abs{\sum_{k=1}^n z_k^\nu} =\Oh{\sqrt n}.
\end{gather}
Instead we proved the related result
\begin{gather} \label{qef}
  \inf_{\abs{z_k} \geq 1}
  \max_{\nu=1,\ldots,\lfloor n^2-n^{65/42+\epsilon} \rfloor} \abs{\sum_{k=1}^n
    z_k^\nu} =\sqrt n + \Oh{n^{23/84+\epsilon}}. \qquad (\epsilon>0)
\end{gather}
In this paper we will improve on these estimates for a general integer $n$.
First we will present a new variant of a construction of Hugh Montgomery (Montgomery's classical construction is given in Tur\'an
\cite{Turan} page 83, or Montgomery \cite{Montgomery}, page 101) to show that \eqref{ajk} is in fact true. We prove that
\begin{gather}
      \inf_{\abs{z_k} \geq 1} \max_{\nu=1,\ldots,\lfloor (m-\varepsilon(n)) n^2 \rfloor } \abs{\sum_{k=1}^n
    z_k^\nu} \leq \sqrt {(m+\varepsilon(n))n},
\end{gather}
for some $\varepsilon(n)=o(1)$. We then use a probabilistic method to
obtain sharper results for Tur\'an's problem 10 proper. We prove the
following theorem.
\begin{thm} \label{ett} One has that
  \begin{gather*}
    \inf_{\abs{z_k} \geq 1}  \max_{\nu=1,\ldots,n^2} \abs{\sum_{k=1}^n z_k^\nu}
    =\sqrt{n}+\Oh{n^{0.2625+\epsilon}}. \qquad \qquad(\epsilon>0)
  \end{gather*}
\end{thm}

\section{Lower bounds}

By using a clever argument involving the Newton-Girard identities,
Cassels \cite{Cassels} proved the following result on pure power sums.
\begin{lem}\label{Cassels} {\rm (Cassels)} Let $(z_1,\ldots,z_n)$ be an $n$-tuple of complex
  numbers. Then
  \begin{gather*}
    \max_{\nu=1,\dots,2n+1} \Re \left( \sum_{k=1}^n z_k^\nu \right)
    \geq 0.
  \end{gather*}
\end{lem}
\noindent As an application he proved that
\begin{gather}
  \max_{\nu=1,\ldots,2n-1} \frac{\max_\nu \abs{ \sum_{k=1}^n z_k^\nu
    }}{\max_\nu \abs{z_k^\nu}} \geq 1.
\end{gather}
In our paper \cite{Andersson} we used Lemma \ref{Cassels} to prove the following result.
\begin{lem} \label{iii} Let $1 \leq m \leq n $ and $(z_1,\ldots,z_n)$ be an $n$-tuple of complex numbers such that $|z_k| \geq 1$.  Then
  \begin{gather*}
    \max_{\nu=1,\dots,2nm-m(m+1)+1} \left| \sum_{k=1}^n z_k^\nu
    \right| \geq \sqrt m.
  \end{gather*}
\end{lem}
\noindent For the special case $n=m$ this implies
\begin{gather} \label{pos}
  \max_{\nu=1,\dots,n^2-n+1} \left| \sum_{k=1}^n z_k^\nu \right| \geq
  \sqrt n, 
\end{gather} and that the lower bound in \eqref{ii} is
valid for all integers $n$. We remark here that there also exist
another theorem independently proved by Newman, Cassels and Szalay
\cite{Szalay} (see Theorem 7.3 in Tur\'an \cite{Turan}), which
furthermore assumes $|z_k|=1$ which in the pure power sum case reduces
to
\begin{lem} {\rm{(Newman, Cassels, Szalay)}} \label{tva} Let $m \geq n$ be an integer. One then has that
\begin{gather*} 
  \max_{1 \leq \nu \leq m} \left| \sum_{k=1}^n z_k^\nu \right| \geq
  \sqrt{n \left(1-\frac {n-1}{m} \right)}.
\end{gather*}
\end{lem}
Although this gives much better estimates than Lemma \ref{iii} for $m
\sim \alpha n^2$ and $0<\alpha<1$, it will give slightly worse estimates
for $m=n^2$. In fact it will give the same estimate when $m \to \infty$
as Lemma \ref{iii} gives for $m=n^2$.

\section{Arithmetical constructions and upper bounds} \label{tre}

In a recent paper \cite{Andersson2} we showed that the three
constructions of Montgomery (\cite{Turan} p. 83), Andersson
\cite{Andersson2} and Fabrykowski \cite{Fabrykowski} gives us the
following estimates.
\begin{lem} \label{expl} One has that 
  \begin{enumerate}[(i)]
  \item if $n+1$ is a prime then there exist an $n$-tuple $(z_1,\ldots,z_n)$ of unimodular complex numbers such that $\abs{s_\nu}\leq \sqrt{n+1} $ for
    $\nu=1,\ldots,n^2+n-1$.
  \item if $n$ is a prime power then there exist an $n$-tuple $(z_1,\ldots,z_n)$ of unimodular complex numbers such that $\abs{s_\nu}\leq \sqrt{n} $ for
    $\nu=1,\ldots,n^2-2$.
  \item if $n-1$ is a prime power then there exist an $n$-tuple $(z_1,\ldots,z_n)$ of unimodular complex numbers such that  $\abs{s_\nu}\leq \sqrt{n-1} $
    for $\nu=1,\ldots,n^2-n$.
  \end{enumerate}
\end{lem}
The main result of our paper \cite{Andersson2} was that Lemma \ref{expl}  $(ii)$ together with Lemma \ref{iii} allows us to obtain explicit solutions to the inf max problem when $n$ is a prime and we take the maximum over $\nu=1,\ldots,n^2-2$ instead of $\nu=1,\ldots,n^2$.

These constructions will also  give us upper estimates in Tur\'an's
problem 10 for all integers, although the inequality will in general not be as sharp. Montgomery's construction, Lemma \ref{expl} $(i)$ gives us slightly  sharper estimates than Lemma \ref{expl} $(ii)$ and $(iii)$.
The construction of Montgomery depends on elementary number theory (
Gauss sums ) whereas the constructions of Fabrykowski
\cite{Fabrykowski} (Lemma \ref{expl} $(iii)$) and Andersson \cite{Andersson2}  (Lemma \ref{expl} $(ii)$) depends on theorems
of Singer \cite{Singer} and Bose \cite{Bose} on projective and affine
geometry over finite fields. We also remark here that all three
constructions are also used in the theory of Sidon sets. For further discussion,  see Andersson \cite{Andersson2} page 6 or Montgomery \cite{Montgomery} page 105-106,  and Martin-O'Bryant
\cite{Martin} for the theory of Sidon sets.

\section{The case of a general integer $n$}

To our knowledge there exist no construction as in Section \ref{tre} for a
general integer $n$. Hence we will try to reduce the general case $n$
to the prime case $p$ by using prime density estimates. We will
consider two variants of this method
\begin{enumerate}
\item Choose a $p \leq n$. Choose $z_1,\ldots,z_p$ by the construction
  given in Lemma 1 and $z_{p+1}, \ldots, z_n$ by the Erd\H os-Renyi
  Lemma.
\item Choose a $p \geq n$. Choose $z_1,\ldots,z_p$ by the construction
  given in Lemma 1. Choose a subset $z_{p_{k_1}}, \ldots,
  z_{p_{k_n}}$.
\end{enumerate}

The difficult part in case 2 will be to estimate the relevant power sums. 
In Section \ref{har} we use trivial methods ( the triangle inequality ). In this case we will only obtain good estimates for a integer $n$ under some strong conjecture on the distribution of primes such as the Cram\'er conjecture. Hence to prove Theorem 1 we need deeper methods and we will use probabilistic methods and moments to show that the desired properties will be true for a random subset.

\subsection{A modified Tur\'an problem 10}
For a general integer $n$ choose the first prime $n \leq p$. Use
Montgomery's construction Lemma \ref{expl} $(i)$
 (Alternatively we can use Lemma \ref{expl} $(ii)$ or $(iii)$) to find an $n$-tuple. Then use the Erd{\H o}s-Renyi Lemma (Lemma \ref{erd}) and choose an $n-p$ tuple $z_{p+1},
\ldots, z_n$ such that
\begin{gather}
  \max_{\nu=1,\dots,p^2} \abs{\sum_{k=p+1}^{n} z_k^\nu} \leq \sqrt {6
    (n-p) \log (p^2+1)}.
\end{gather}
By the triangle inequality it is clear that
\begin{align*}
  \max_{\nu=1,\dots,p^2} \abs{\sum_{k=1}^{n} z_k^\nu} & \leq
  \max_{\nu=1,\dots,p^2} \abs{\sum_{k=1}^{p} z_k^\nu} +
  \max_{\nu=1,\dots,p^2} \abs{\sum_{k=p+1}^{n} z_k^\nu}, \\ &\leq
  \sqrt{p+1} + \sqrt {6 (n-p) \log (p+1)}.
\end{align*}
This idea was first used by Queffelec \cite{Queffelec} and yields the upper bound in the following Proposition. The lower bound can be obtained similarly with equation \eqref{pos}.
\begin{prop}
  Let $p$ be a prime and $n \geq p$ an integer. Then
  \begin{gather*}
    \sqrt{p} - \sqrt {6 (n-p) \log (p+1)}  \leq \inf_{|z_k|=1}\max_{\nu=1,\dots,p^2} \abs{\sum_{k=1}^{n} z_k^\nu}
    \leq \sqrt{p+1} + \sqrt {6 (n-p) \log (p+1)}.
  \end{gather*}
\end{prop}

By combining this with the estimate of Baker-Harman-Pintz
\cite{BakerHarman} on the difference between consecutive primes
\begin{lem} \label{bhp}
  {\rm (Baker-Harman-Pintz)} Suppose that $p_k$ denote the $k$'th prime.
  Then
  \[ p_{k+1}-p_k \leq p_k^{0.525}. \qquad \qquad (p_k>N_0) \]
\end{lem}

\noindent we obtain the following Proposition.
\begin{prop} \label{prop1} One has that
  \begin{gather*}
    \inf_{\abs{z_k} = 1} \max_{\nu=1,\dots,n^2-2n^{1.525}}
    \abs{\sum_{k=1}^{n} z_k^\nu} = \sqrt n + \Oh{n^{0.2625+\epsilon}}. \qquad \qquad (\epsilon>0)
  \end{gather*}
\end{prop}
We remark here that the reason why we get $0.2625$ instead of $23/84=0.273809\ldots$ as we had in \cite{Andersson} is that we use sharper estimates of the difference between consecutive primes, Baker-Harman-Pintz \cite{BakerHarman} instead of the estimate of Iwaniec-Pintz \cite{IwaniecPintz}.  

A problem with Lemma \ref{bhp} is that it will not give us an estimate for
Tur\'an's problem 10. In order to get such an estimate it is sufficient
to have an explicit construction that allows us to choose the maximum
over $\nu=1,\ldots,\lfloor \alpha n^2 \rfloor$ for some $\alpha>1$. In subsection \ref{montg}  we will see how a new variant of Montgomery's construction
will allow us to choose  any $\alpha>1$.

\subsection{A problem from operator theory}

We  remark here that  Proposition \ref{prop1} also has an application on operator theory. Let $k_n$ be the smallest constant such that for any $n$-dimensional normed space $X$ and any invertible operator $T \in \mathcal  L (X)$ we have that
 \begin{gather*}
     \abs{\det(T)} \norm{T^{-1}} \leq k_n \norm{T}^{n-1}
  \end{gather*}
Sch\"affer \cite{Schaffer} proved that $k_n \leq  \sqrt {en}$. Gluskin-Meyer-Pajor \cite{GluskinMeyerPajor} who seemed unaware of Tur\'an's book \cite{Turan} and Erd{\H o}s-Renyi's paper \cite{Erdos} obtained an independent proof of Erd{\H o}s-Renyi's result Lemma \ref{erdren} and used it to prove that $k_n \geq c \sqrt {n/\log n}$. Queffelec \cite{Queffelec}  used  Gluskin-Meyer-Pajor's method but substituted the use of an Erd\H os-Renyi type result  with a variant of Proposition \ref{prop1} ( equation \eqref{qef} ) to prove that $k_n \geq \sqrt n (1-o(1))$.
  He did not either refer to Tur\'an's book and obtained Montgomery's construction independently. For further results on this problem, see Nikolski \cite{Nikolski}. Since we have that
\begin{gather*}
  \sqrt n (1-o(1)) \leq  k_n \leq  \sqrt {en}
\end{gather*}
it seems reasonable to state the following conjecture.
\begin{conj} \label{operator}
  There exist a constant $1 \leq c \leq \sqrt e$ such that
 \begin{gather*}
  k_n \sim c \sqrt n.
\end{gather*}
\end{conj}
\begin{prob} Solve  Conjecture \ref{operator} and find the constant $c$. \end{prob}
By following the proof of Theorem 4 in Gluskin-Meyer-Pajor \cite{GluskinMeyerPajor} p. 235 (which they attribute to J. Bourgain) it does not seem as if known results from the Tur\'an power sum method can improve on the lower bound which is essentially $\sqrt n$.

 The proof of the upper bound uses completely different methods (operator theoretic), that does not seem easy to sharpen  as well. The key point in a possible proof of Conjecture \ref{operator} might very well be identity (1) in \cite{GluskinMeyerPajor}.

By studying the proof of Theorem 4  \cite{GluskinMeyerPajor} p. 235  more carefully it is clear  that if for each $\epsilon>0$ and some sufficiently large $n \geq n_0(\epsilon)$ there exist an $n$-tuple $(z_1,\ldots,z_n)$ of unimodular complex numbers and a $c$ such that
 \begin{gather*}
     \max_{\nu=1,\dots,\lfloor n \log n \rfloor}   \abs{\sum_{k=1}^{n} z_k^\nu} \leq c \sqrt n (1+\epsilon), 
  \end{gather*}
then 
\begin{gather*}
  k_n \geq \frac {\sqrt n} c.
\end{gather*}
By Sh\"affer's result $k_n \leq \sqrt{en}$ we find that $c \geq 1/\sqrt e$. Hence we see that Sh\"affer's result combined with the theory from Gluskin-Meyer-Pajor \cite{GluskinMeyerPajor}, in particular the proof of Bourgain will give a third method to prove lower bounds in Tur\'an's problem. One obtains that
  \begin{gather} \label{r}
    \sqrt{\frac n e}(1-o(1))  \leq  \inf_{\abs{z_k} = 1} \max_{\nu=1,\dots,\lfloor n \log n \rfloor}   \abs{\sum_{k=1}^{n} z_k^\nu}.   
  \end{gather}
We remark that this gives worse lower estimates than Lemma \ref{tva} which
implies that $\sqrt {n/e}$ can be replaced by $\sqrt n$ in equation \eqref{r}. The proof of Lemma \ref{tva} is much simpler as well, nevertheless we find it amusing that results from operator theory implies results in Tur\'an power sum theory.

\subsection{Montgomery's construction}\label{montg} 
 
Let $\chi$ be a character mod $p$, and $\chi_0$ the trivial character.
From the theory of Gauss sums we have that
\begin{gather} \label{aaa} \abs{\sum_{k=1}^{p-1} \chi(k) e
    \p{\frac{ka} p}} =
  \begin{cases} \sqrt{p}, & p \not | a, \, \chi \neq \chi_0, \\ 1, & p
    \not |a, \, \chi=\chi_0, \\ 0, & p|a, \, \chi \neq \chi_0, \\ p-1, &
    p| a, \, \chi=\chi_0. \end{cases}
\end{gather}
The idea of Montgomery (see Tur\'an
\cite{Turan} page 83 or Montgomery \cite{Montgomery}, page 101) is to use
\begin{gather} \label{ur1} z_k=\chi(k) e \p{\frac{k} p}, \qquad \qquad
  (k=1,\ldots,p-1)
\end{gather}
where $\chi$ is a character mod $p$ of order $p-1$. Lemma \ref{expl} $(i)$
now follows from \eqref{aaa}.

We now assume that $p=nm+1$ and $\chi$ is  character mod
$p$ of order $p-1$, and let
\begin{gather}
  w_k=\chi(k) e \p{\frac{k} {nm+1}}, \qquad \qquad (k=1,\ldots,mn) \\
  \intertext{and}
  \label{ur12}
  z_k=w_k^m. \qquad \qquad (k=1,\ldots,n)
\end{gather}
In other words this means that $\{z_k \}$ is the subset of $\{w_j\}$
where $j$ is an $m$'th power residues of $p$. It is clear that
\begin{gather*}
  \frac 1 m \sum_{j=1}^{m} \chi(k)^{nj} = \begin{cases} 1, & k \text{
      is a $m$'th power residue,} \\ 0, & \text{otherwise.} \end{cases}
\end{gather*}
Hence
\begin{gather*} \begin{split}
    \sum_{k=1}^n z_k^\nu &= \sum_{k=1}^{nm} \frac 1 m \sum_{j=1}^{m} \chi(k)^{nj} w_{k}^{\nu}, \\
    &= \frac 1 m \sum_{j=1}^{m} \sum_{k=1}^{nm} \chi(k)^{nj+\nu} e
    \p{\frac{k \nu} {nm+1}}.
  \end{split}
\end{gather*}
By \eqref{aaa} each term except when $nm | nj+\nu$ for $j=1,\ldots,m$
and $(nm+1)|\nu$ will contribute at most $\sqrt{nm+1}$. Since $nm |
nj+\nu$ implies that $n| \nu$, and $p| \nu$ furthermore implies that
$pn | \nu$ this implies that unless $(mn+1)n | \nu$ we have that
\begin{gather*}
  \abs {\sum_{k=1}^n z_k^\nu} \leq \frac 1 m \sum_{j=1}^{m}
  \abs{\sum_{k=1}^{nm} \chi(k)^{nj+\nu} e \p{\frac{k \nu} {nm+1}}}
  \leq \sqrt{nm+1}.
\end{gather*}
Hence the following Proposition follows.
\begin{prop} \label{montmod}
  Suppose that $p=mn+1$ is prime. Then there exist an $n$-tuple $(z_1,\ldots,z_n)$ of unimodular complex
  numbers such that
  \begin{gather*}
    \abs{\sum_{k=1}^n z_k^\nu} \leq \sqrt{mn+1}. \qquad
    (\nu=1,\ldots,mn^2+n-1)
  \end{gather*}
\end{prop}
By the same method as used to prove Proposition \ref{prop1} we can use
this to prove the proposition.
\begin{prop} Let $m \geq 1$ be an integer and $\epsilon>0$. Then one has that
  \begin{gather*} 
    \sqrt{n}- \Oh{n^{0.275+\epsilon}} \leq \inf_{z_k \in \C, |z_k| = 1} \max_{\nu=1,\dots,
      \lfloor m n (n-n^{0.55+\epsilon}) \rfloor} \left| \sum_{k=1}^n z_k^\nu \right| \leq \sqrt {mn}+ \Oh{n^{0.275+\epsilon}}.
  \end{gather*}
\end{prop}

\begin{rem} \label{harrem}
The reason why we get $0.275$
    instead of $0.2625$ is that we need primes in arithmetical
    progressions, i.e. primes $\equiv 1 \pmod m$, and instead of the
    Baker-Harman-Pintz theorem \cite{BakerHarman} we can use the
    Baker-Harman-Pintz theorem \cite{BakerHarmanII} for arithmetical
    progressions.  Since each odd prime $\equiv 1 \pmod 2$ this means
    that in the special case $m=2$ we can instead use $0.2625$. Professor Glyn Harman has informed us that by using a (although not very effective in $m$) method from Harman-Watt-Wong \cite{HarmanWattWong} they can obtain the same constant $0.525$ for the difference between consecutive primes in an arithmetical progression as for general primes. Hence the constant $0.275$ in Proposition 4 can be replaced by $0.2625$.
\end{rem}

In particular this implies that for the Tur{\'a}n problem (by also using the lower bound \eqref{pos}) we have that
\begin{gather*}
  \sqrt{n} \leq \inf_{\abs{z_k} \geq 1} \max_{v=1,\ldots,n^2}\abs{\sum_{k=1}^n z_k^v} \leq
  \sqrt{2 n}(1+ o(1)),
\end{gather*}
which improves on the bound of Erd\H os and Renyi, equation \eqref{erdren}.

\section{A general integer II} \label{har}

\begin{prop} \label{prop10} One has that 
\begin{enumerate}[(i)] \item  Suppose that $p=n+j+1$ is a prime for $j \geq 0$. Then

  $\qquad \qquad \displaystyle \sqrt n \leq \inf_{\abs{z_k} \geq 1} \max_{\nu=1,\ldots,n^2}
    \abs{\sum_{k=1}^n z_k^\nu}\leq \sqrt{n+j+1} +j,$ \item Suppose that
      $n+j$ is a prime power for $j \geq 0$. Then  

$\qquad \qquad \displaystyle \sqrt n \leq \inf_{\abs{z_k}
      \geq 1} \max_{\nu=1,\ldots,n^2} \abs{\sum_{k=1}^n z_k^\nu}\leq
    \sqrt{n+j} + j,$  \item  Suppose that $n$ is a prime power.  Then

  $\qquad \qquad \displaystyle \sqrt n \leq \inf_{\abs{z_k} \geq 1} \max_{\nu=1,\ldots,n^2}
    \abs{\sum_{k=1}^n z_k^\nu}\leq \sqrt{n} + 1.ö$
  \end{enumerate}
\end{prop}
\begin{proof}
We first prove $(i)$.  The lower bound follows from \eqref{pos}. The upper bound follows from the following construction. Choose $z_1,\ldots,z_{n+j}$ by the
construction given in Lemma \ref{ett} $(i)$. By the triangle inequality it is clear that
\begin{gather*}\begin{split}
    \abs{\sum_{k=1}^{n} z_k^\nu} &= \abs{\sum_{k=1}^{n+j}
      z_k^v-\sum_{k=n+1}^{n+j} z_k^\nu }, \\ &\le \abs{\sum_{k=1}^{n+j}
      z_k^\nu} + \abs{\sum_{k=n+1}^{n+j} z_k^\nu}, \\ &\le \sqrt{n+j+1}+ j.
  \end{split}
\end{gather*}
  The proof of $(ii)$ and $(iii)$ follows by using Lemma \ref{ett} $(ii)$ and $(iii)$ instead of  Lemma \ref{ett} $(i)$.
\end{proof}
\begin{rem}
 The reason why we have stated Proposition \ref{prop10} $(iii)$ for $n+j-1$ prime power only in the case $j=1$ is that for all $j \geq 1$ Proposition \ref{prop10} $(ii)$ will give sharper results. Similarly it is easily seen that if $n+j$ in  Proposition \ref{prop10} $(ii)$ is in addition to being a prime power also prime, then  Proposition \ref{prop10} $(i)$ will give sharper results.
\end{rem}

By the Cram\'er conjecture \cite{Cramer}
\begin{conj} {\rm (Cram\'er)} $p_{k+1}-p_k = \Oh{ (\log p_k)^2}$
\end{conj}
\noindent we obtain the following conditional result.
\begin{prop}
  The Cram\'er conjecture implies that
  \begin{gather*}
    \inf_{\abs{z_k} \geq 1} \max_{\nu=1,\ldots,n^2}
    \abs{\sum_{k=1}^n z_k^\nu}= \sqrt{n} + \Oh{ (\log n)^2}.
  \end{gather*}
\end{prop}
By the Riemann hypothesis it follows that (Cram\'er \cite{Cramer})
\begin{gather*}
  p_{k+1}-p_k \ll \sqrt p_k \log p_k,
\end{gather*}
and we see that even under the Riemann hypothesis, Proposition \ref{prop10} does not give any better result than Erd{\H o}s-Renyi's  result equation \eqref{erdren}. Hence if we like to use Proposition \ref{prop10} to prove asymptotic estimates in Tur\'an's problem 10 we need a stronger estimate for the distribution of consecutive primes, such as the Cram\'er conjecture, or at least $p_{k+1}-p_k=\Oh{p_k^\theta}$ for some $\theta<1/2$. Since no such result exists we will seek other methods of proof.

\section{Moments}
\subsection{A fundamental lemma}
We will first prove a more technical lemma before we prove our main lemma. Let $(z_1,\ldots,z_n)$ be an $n$-tuple of complex numbers. Define
\begin{gather} \label{stst2}
  S(\nu_1,\ldots,\nu_m)=\sum_{\substack{k_1,\ldots,k_m=1 \\ i \neq j \implies k_j \neq k_i}}^n z_{k_1}^{\nu_1} \cdots  z_{k_m}^{\nu_m}.
\end{gather}
We see that for $m=1$ this reduces to the classical power sum
\begin{gather} \notag
  S(\nu)= \sum_{k=1}^k z_k^\nu.
\intertext{It is clear that} \label{stst}
S(\nu_1) \cdots S(\nu_m)=\sum_{k_1,\ldots,k_m=1}^n z_{k_1}^{\nu_1} \cdots  z_{k_m}^{\nu_m}.
\end{gather}
We recall that  $U$ is a disjoint union of nonempty sets of $S$ if $U$ is a family of sets $\{U_i:i=1,\ldots,k \}$ such that $U_i \neq \emptyset, S=\cup_{i=1}^k U_i$ and $U_i \cap U_j=\emptyset$ for $i \neq j$. We can expand \eqref{stst} in terms of  \eqref{stst2} as follows.
\begin{lem} \label{lemma7}
  Let $\mathcal U_m=\{U=\{U_1,\ldots,U_k\} \} $  be the family of all disjoint unions of nonempty sets of  $\{1,\ldots,m\}$.  Then one has that
   \begin{gather*}
    S(\nu_1) \cdots S(\nu_m)= \sum_{U \in \mathcal U_m}  S\left(\sum_{j \in U_1} \nu_j,\ldots,\sum_{j \in U_k} \nu_j \right).
 \end{gather*}
\end{lem}
 Furthermore one has that the elements   of $S(\nu_1) \cdots S(\nu_m)$ where the product contains  exactly $k$ different elements $z_{j_1},\cdots,z_{j_k}$ can be written as

  \begin{gather} \label{trt}
    \sum_{\substack{U \in \mathcal U_{m} \\ \# U=k}}  S\left(\sum_{j \in U_1} \nu_j,\ldots,\sum_{j \in U_k} \nu_j \right).
  \end{gather}

We have the following lemma.
\begin{lem} \label{prefund} Let $(z_1,\ldots,z_n)$ be an $n$-tuple of unimodular complex numbers.  Suppose that $\sqrt n \leq M(k)$ and $\abs{S(k)}  \leq M(k)$ for all integers $k$. Then
  \begin{gather*}
    \abs{S(\nu_1,\ldots,\nu_m)} \ll C_m M(\nu_1) \cdots M(\nu_m).
  \end{gather*}
\end{lem}
\begin{proof} 
  We will use the principle of induction. By the assumption $\abs{S(k)} \leq M(k)$ the assertion is true for $m=1$. Now assume it is true for $m=m_0$.  Lemma \ref{lemma7} gives us
  \begin{gather}
    S(\nu_1) \cdots S(\nu_{m_0+1})= \sum_{U \in \mathcal U_{m_0+1}}  S \left(\sum_{j \in U_1} \nu_j,\ldots,\sum_{j \in U_k} \nu_j \right)
 \end{gather}
 There is a unique element in $\mathcal U_{m_0+1}$ with $m_0+1$ elements $\{\{1\},\ldots,\{m_0+1\} \}$. We see that
 \begin{gather*}
   S(\nu_1,\ldots,\nu_{m_0+1})  = S(\nu_1) \cdots S(\nu_{m_0+1}) -
  \sum_{\substack{U \in \mathcal U_{m_0+1} \\ \# U  \leq m_0}} S \left(\sum_{j \in U_1} \nu_j,\ldots,\sum_{j \in U_k} \nu_j \right)
 \end{gather*}
 The first part is less than  $M(\nu_1) \cdots M(\nu_{m_0+1})$ by the assumption $\abs{S(\nu_j)} \leq M_{\nu_j}$. By the assumption the Lemma is true for $m=m_0$ and the sum over $\mathcal U_{m_0+1}$ over disjoint unions of $S$ with at most $m_0$ elements can by be  estimated by
\begin{gather*}
   \sum_{\substack{U \in \mathcal U_{m_0+1}\\ \# U  \leq m_0}} M\left(\sum_{j \in U_1} \nu_j\right) \cdots M\left(\sum_{j \in U_k} \nu_j \right)  
\end{gather*}
By the argument $\sqrt n \leq M(k)$ and the trivial fact that $\abs{S(k)}\leq n$, it follows that this as well can be estimated by $C M(\nu_1) \cdots M(\nu_{m_0+1})$. Hence the Lemma is true for $m_0+1$. The general results follows from the principle of induction.
\end{proof}

\begin{lem} \label{fund} {\rm (Fundamental Lemma)}  Let $\alpha,\epsilon>0$, $0<\theta<1$ and $C_1 \geq 1$ be given. Suppose that  $(z_1,\ldots,z_n)$ is an $n-$tuple of  unimodular complex numbers,
 \begin{gather} \label{star2}
  m \sim n^\theta, 
 \end{gather}
the quantity $S(\nu)$ denote the pure power sum
  \begin{gather*}
    S(\nu)=\sum_{k=1}^n z_k^\nu, \\ \intertext{and}
    \abs{S(\nu)} \leq C_1 \sqrt n. \qquad \qquad (\nu=1,\ldots, \lfloor \alpha n^2 \rfloor) \\ \intertext{Let $\mathcal N=\{1,\ldots,n \}$. Then there exist a subset $\mathcal M_0  \subset  \mathcal N$, with $\# \mathcal M_0=m$  such that 
 }
    \abs{\sum_{k \in \mathcal M_0} z_{k}^\nu} \ll_{\epsilon} m^{1/2+\epsilon}. \qquad  \qquad(\nu=1,\ldots, \lfloor \alpha n^2 \rfloor)
  \end{gather*}
\end{lem}
\begin{proof}

In order to find  the subset $\mathcal M_0 \subset \mathcal M$  we  use probabilistic methods (moments). We choose an integer 
\begin{gather}\label{star1}
   N > \frac 1 {\theta \epsilon},
\end{gather}
and consider the sum over all subsets
\begin{gather*}
    \sum_{\substack{\mathcal M \subset \mathcal N \\ \# \mathcal M =m}}\abs{\sum_{k \in \mathcal M} z_k}^{2N}. \qquad \qquad (\mu=1,\ldots,\lfloor \alpha n^2 \rfloor)
\end{gather*}
There are $\binom{n} m $ such sets. Hence we can choose a subset $\mathcal M_0 \subset \mathcal N$ with $\# \mathcal M_0=m$  such that
\begin{gather} \label{jag}
    \frac 1 {\binom n m} \abs{\sum_{k \in \mathcal M_0} z_k^\nu}^{2N} \leq  \sum_{\mu=1}^{ \alpha n^2} \sum_{\substack{\mathcal M \subset \mathcal N \\ \# \mathcal M=m}} \abs{\sum_{k \in \mathcal M} z_k^\mu}^{2N}
\end{gather}
for each $\nu=1,\ldots,\lfloor \alpha n^2 \rfloor$. We consider
\begin{gather} \label{ajaj}
  \frac 1 {\binom n m} \sum_{\substack{\mathcal M \subset \mathcal N \\ \# \mathcal M =m}} \abs{\sum_{k \in \mathcal M} z_k^\mu}^{2N}  =   
  \frac 1 {\binom n m} \sum_{\substack{\mathcal M \subset \mathcal N \\ \# \mathcal M =m}}  \sum_{k_1,\ldots,k_{2N} \in \mathcal M} 
   (z_{k_1} \cdots z_{k_N})^\mu  (z_{k_{N+1}} \cdots z_{k_{2N}})^{-\mu}.
 \end{gather}

As the sum is over $\mathcal M  \subset \mathcal N$  each term can be written as
\begin{gather} \label{aj}
    (z_{k_1} \cdots z_{k_n})^\mu \p{z_{k_{N+1}} \cdots z_{k_{2N}}}^{-\mu}
\end{gather}
with $k_i \subset \mathcal N$.  Suppose that $\mathcal K= \{k_1,\ldots, k_{2N} \}$ and $k=\# \mathcal K$. By a simple combinatorial argument  we can choose
 \begin{gather*}
   \binom{n} {m-k}
 \end{gather*}
subsets $\mathcal M \subset \mathcal N$ such that $\mathcal K \subset \mathcal M$. Hence each term \eqref{aj} will occur with the coefficient $\binom n {m-k}$ and \eqref{ajaj}
 can be written as
 \begin{gather}  \label{utr}
   \frac 1 {\binom n m}  \sum_{k_1,\dots,k_{2N} \in \mathcal N} \binom{n}{m-\#\{k_1,\dots,k_{2N} \}} 
   \p{z_{k_1} \cdots z_{k_N}}^{\mu} \p{z_{k_{N+1}} \cdots z_{k_{2N}}}^{-\mu}.
\end{gather}
By equation \eqref{trt}  this equals
\begin{gather*}
    \sum_{k=1}^{2N} \frac{\binom{n}{m-k }} {{\binom n m}} 
\sum_{\substack{U \in {\mathcal U}_{2N} \\ \# U=k}}  S\left( \sigma(U_1) \mu ,\ldots, \sigma(U_k) \mu \right),
\end{gather*}
where  $\sigma(S)=\# S \cap \{1,\ldots,N \} -\# S \cap \{N+1,\ldots,2N \}$, $S(\nu_1,\ldots,\nu_{2N})$ is defined by eq. \eqref{stst2} and $\mathcal U_{2N}$ is defined in Lemma \ref{lemma7}. By Lemma \ref{prefund} this can be estimated by
\begin{gather*}
    \sum_{k=1}^{2N} \frac{\binom{n}{m-k }}   {{\binom n m}}\sum_{\substack{U \in {\mathcal U}_{2N} \\ \# U=k}} M(\sigma(U_1) \mu)  \cdots M  ( \sigma(U_k)\mu).
\end{gather*}
If $U_j$ has $1$ element than $M( \sigma(U_k) \mu) = M(\pm \mu)=C_1 \sqrt n$. If $ \# U_j \geq 2$, then $M(\sigma(U_j) \mu) \leq n$.
This implies that
\begin{gather*}
 \sum_{\substack{U \in {\mathcal U}_{2N} \\ \# U=k}} M(\sigma(U_1) \mu)  \cdots M  ( \sigma(U_k)\mu) \ll n^{\min(N,k)}.
\end{gather*}
By the further fact that
\begin{gather*}  \frac {\binom{n}{m-k}}  {\binom n m } \leq   \p{\frac m n}^k,
\\ \intertext{this implies  that equation \eqref{utr} can be estimated by}
    \sum_{k=1}^{2N} \p{\frac m n}^k n^{\min(N,k)}. 
\end{gather*}
The dominating term will be $k=N$ and this can be estimated by 
\begin{gather*}
  \Oh{m^N}. \qquad \qquad (\mu=1,\ldots,\lfloor \alpha n^2 \rfloor)
\end{gather*}
When we sum over $\mu=1,\ldots \alpha n^2$  in \eqref{jag} we get that
\begin{gather*}
  \abs{\sum_{k \in \mathcal M_0} z_k^\nu}^{2N} \ll \alpha n^{2} m^N.  \qquad   \qquad (\nu=1,\ldots,\lfloor \alpha n^2 \rfloor) 
\end{gather*}
By equations \eqref{star2} and  \eqref{star1} we obtain
\begin{gather*}
  \abs{\sum_{k \in \mathcal M_0} z_k^\nu}  \ll m^{1/2+\epsilon}.  \qquad \qquad (\nu=1,\ldots,\lfloor \alpha n^2 \rfloor) 
\end{gather*}
\end{proof}

\subsection{Proof of Theorem 1} 
\begin{proof}
The lower bound follows from equation \eqref{pos}. Hence we will concentrate on the upper bound. By the Baker-Harman-Pintz theorem, Lemma \ref{bhp} we can choose a prime $n<p$ such that $p-n \asymp p^{0.525}$. 
By the Montgomery construction, Lemma \ref{expl} $(i)$ (Or alternatively, we can use Lemma \ref{expl} $(ii)$ or $(iii)$). we can choose a $(p-1)-$tuple $(z_1,\ldots,z_{p-1})$   of unimodular complex numbers such that
\begin{gather*}
  \abs{\sum_{k=1}^{p-1} z_k^\nu} \leq  \sqrt {p}. \qquad  \qquad (\nu=1,\ldots,(p-1)^2) 
\end{gather*}
Let $m=p-1-n$. By the fundamental Lemma \ref{fund} with $\alpha=1$ and $\theta=0.525$ we can choose a subset $\mathcal M_0 \subset \{1,\ldots,p-1\}$ with $\# M_0=m$ such that
\begin{gather*}
  \abs{\sum_{k \in \mathcal M_0} z_k^\nu} \leq  n^{0.2625+\epsilon}. \qquad \qquad (\nu=1,\ldots,(p-1)^2) \qquad (\epsilon>0)
\end{gather*}
Let $\mathcal N = \{1,\ldots,p-1 \} \setminus \mathcal M_0$. It is clear that $\# \mathcal N=n$ and by the triangle inequality it follows for $1 \leq \nu \leq n^2 \leq (p-1)^2$ that
\begin{gather*}
     \begin{split}
     \abs{\sum_{k \in \mathcal N} z_k^\nu} &=  \abs{\sum_{k=1}^{p-1} z_k^\nu - \sum_{k \in \mathcal M_0} z_k^\nu}, \\ 
        &=  \abs{\sum_{k=1}^{p-1} z_k^\nu}+ \Oh{\abs{\sum_{k \in \mathcal M_0} z_k^\nu}},  \\
           &=  \sqrt{n+\Oh{n^{0.525}}} + \Oh{n^{0.2625+\epsilon}}, \\ &= 
             \sqrt n + \Oh{n^{0.2625+\epsilon}}.
   \end{split}
\end{gather*}
which finishes the proof of Theorem 1. \end{proof}

By the same proof method but with the  modified Montgomery construction Proposition \ref{montmod} instead of the
 classical Montgomery construction, and the Baker-Harman-Pintz theorem for primes in arithmetical progressions \cite{BakerHarmanII} we obtain the following result.
\begin{thm} \label{gmt}
  Let $m \geq 1$ be an integer.  One then has that
  \begin{gather*} 
    \sqrt{n} \leq \inf_{z_k \in \C, |z_k| \geq 1} \max_{\nu=1,\dots,
       m n^2} \left| \sum_{k=1}^n
      z_k^\nu \right| \leq \sqrt {mn}+ \Oh{n^{0.275+\epsilon}}. \qquad (\epsilon>0)
  \end{gather*}
\end{thm}
\begin{rem}
 As in remark \ref{harrem} the constant $0.275$ in Theorem \ref{gmt} be replaced by $0.2625$
\end{rem}
\section{Tur\'an's problem 10 on the average}

In Section \ref{har}  we proved conditional results 
(under the Cram\'er conjecture). In this section we will show sharper results on the average.  Let $\Delta(n)$ be defined by
\begin{gather} \label{Deltadef}
    \inf_{z_k \in \C, |z_k| \geq 1} \max_{\nu=1,\dots, n^2} \left| \sum_{k=1}^n z_k^\nu \right|=\sqrt n+ \Delta(n).
\end{gather}
Theorem \ref{ett} and the positivity eq. \ref{pos} gives us
\begin{gather} 0 \leq \Delta(n) \ll n^{0.2625+\epsilon}. \qquad \qquad (\epsilon>0)
\end{gather}
In more generality we have that the proof method of Theorem \ref{ett} and Lemma \ref{fund} implies that if
$p_k \leq n \leq p_{k+1}$ for consecutive primes, then we have that
\begin{gather} \label{ut9}
  \Delta(n) \ll n^{\epsilon} \sqrt{p_{k+1}-p_k}. \qquad \qquad (\epsilon>0)
\end{gather}
From this there follows a number of results on the average order of $\Delta(n)$ by corresponding results for the average orders of differences of consecutive primes. We have the following theorem.
\begin{thm}
 One has that
\begin{gather*}
  \sum_{n=1}^N \abs{\Delta(n)}^2 \ll N^{5/4+\epsilon}. \qquad \qquad (\epsilon>0)
\end{gather*}
\end{thm}
\begin{proof} 
This follows from eq. \eqref{ut9} by using the estimate
\begin{gather*}
  \sum_{p_k \leq X} (p_{k+1}-p_k)^2 \ll X^{5/4+\epsilon}\qquad \qquad (\epsilon>0)
\end{gather*}
from Peck's D.Phil Thesis.
\end{proof}
\begin{rem} The constant $5/4$ improves on the constant $23/18$ of Heath-Brown \cite{HeathBrown}. We are grateful to Professor Glyn Harman for informing us of Peck's result. \end{rem}
We can also prove the following result.
\begin{thm} \label{harg}
\begin{enumerate}[(i)]
  \item Under the Density hypothesis one has that that $$\qquad \qquad 0 \leq \Delta(n) \ll n^{1/4+\epsilon}. \qquad \qquad (\epsilon>0)$$
  \item Under the Lindel\"of hypothesis one has that$$  \displaystyle \sum_{n=1}^N \abs{\Delta(n)}^2 \ll N^{1+\epsilon}.\qquad \qquad (\epsilon>0)$$
\end{enumerate}
\end{thm}
\begin{proof}
 The Density hypothesis implies that (see e.g. Ivic \cite{Ivic})
\begin{gather}
   p_{k+1}-p_k \leq  p_k^{1/2+\epsilon}. \qquad \qquad (\epsilon>0) \\
\intertext{Yu \cite{Gang} has proved that  } \sum_{p_k \leq X} (p_{k+1}-p_k)^2 \ll X^{1+\epsilon} 
\end{gather}
 is true under the Lindel\"of hypothesis. Together with equation \eqref{ut9} this implies our Theorem.
\end{proof}

\begin{rem} The Riemann hypothesis implies the Lindel\"of hypothesis, and the Lindel\"of hypothesis implies the Density hypothesis (see e.g. Ivic \cite{Ivic}), hence the statements in Theorem \ref{harg} are true also under the Riemann hypothesis. In this case we could have used the more classical results  that the Riemann hypothesis implies that $p_{k+1}-p_k \leq \sqrt{p_k} \log p_k$  of Cram\'er \cite{Cramer} and  $\sum_{p_k \leq X} (p_{k+1}-p_k)^2 \ll X (\log X)^3$ 
which is a result of Selberg \cite{Selberg} to prove Theorem \ref{harg}
 \end{rem}

 Unconditionally we can use equation \eqref{ut9} and a theorem of Peck \cite{Peck} for how many $k$'s that fulfills $p_{k+1}-p_k \geq \sqrt{p_k}$ to get an estimate for how many $n$'s that does not fulfill this estimate.

\begin{thm} One has that
\begin{gather*}
  \sum_{\substack{1 \leq n \leq X \\ \Delta(n) \gg n^{1/4+\epsilon}}} 1 \ll X^{25/36+\epsilon}.\qquad \qquad (\epsilon>0)
\end{gather*}
\end{thm}

\section{Further problems}

We will here investigate the following problem.
\begin{prob} Let $\alpha>0$ be a real number. Find an asymptotic formula for
\begin{gather*} 
  \inf_{z_k \in \C, |z_k| \geq 1} \max_{\nu=1,\dots, \lfloor
       \alpha n^2 \rfloor} \left| \sum_{k=1}^n z_k^\nu \right|.
 \end{gather*}
\end{prob}
 Theorem \ref{ett} proves strong results for $\alpha=1$. For general values of $\alpha$ the problem of getting true asymptotics seems more difficult.

\begin{thm} \label{u77} One has that
\begin{gather*} 
  (\underline{A} (\alpha)-o(1)) \sqrt n  \leq \inf_{z_k \in \C, |z_k| \geq 1} \max_{\nu=1,\dots,
       \lfloor \alpha n^2 \rfloor} \left| \sum_{k=1}^n z_k^\nu \right| \leq (\overline A(\alpha)+o(1))\sqrt n,
 \\ \intertext{for}
   \underline A(\alpha)= \begin{cases} 1-\sqrt{1-\alpha}, & 0 < \alpha \leq 1, \\ 1, & \alpha>1, \end{cases} \qquad \text{ and } \qquad \overline A(\alpha)= \begin{cases} 1, & 0 <\alpha \leq 1, \\ \sqrt 2, & 1<\alpha \leq 2, \\ \sqrt 3, & 2<\alpha \leq 3, \\ 2, & 3< \alpha. \end{cases} 
\end{gather*}
\end{thm}
\begin{proof}
  The lower bound follows from Lemma \ref{iii}. The upper bound follows from Theorem \ref{gmt} for $1 \leq \alpha \leq 4$. For $\alpha>4$ it follows from a new result of ours, Corollary 2 of our recent paper  \cite{Andersson4}. In fact our paper \cite{Andersson4} answers several questions from version 2 of this paper on arXiv, see \cite{Andersson5}, pages 17-20.
   
\end{proof}

One sees that the only case where we know the true asymptotics is in fact $\alpha=1$, or in other words Problem 10 of Tur\'an which we already studied in more detail.
\begin{conj} \label{con3} One can choose $A(\alpha)=\underline A(\alpha)=\overline A(\alpha)$ in Theorem \ref{u77}. \end{conj}

We tend to believe that $A(\alpha)=1$ for $0<\alpha<1$. The following theorem from our recent paper \cite{Andersson3} proves this under the further assumption that $\abs{z_k}=1$.

\begin{thm} (Andersson, 2007) \label{avs} Let $\alpha >0$ be a constant. One then has that
\begin{gather*}
 (\underline{B} (\alpha)-o(1)) \sqrt n  \leq \inf_{z_k \in \C, |z_k| = 1} \max_{\nu=1,\dots, \lfloor \alpha n^2 \rfloor} \left| \sum_{k=1}^n z_k^\nu \right| \leq (\overline B(\alpha)+o(1))\sqrt n,
\\ \intertext{where}
        \underline{B}(\alpha)= \begin{cases}  1, & 0<\alpha\leq 1, \\ \sqrt{\frac 3 2-\frac 1 {2\alpha}}, & 1 \leq \alpha \leq 3, \\
                                    \sqrt{2-\frac 2 \alpha}, & 3 \leq \alpha. \end{cases} \qquad \text{and} \qquad
 \overline B(\alpha)= \begin{cases} 1, & 0 <\alpha \leq 1, \\ \sqrt 2, & 1<\alpha \leq 2, \\ \sqrt 3, & 2<\alpha \leq 3, \\ 2, & 3< \alpha. \end{cases} 
\end{gather*}
\end{thm}

\begin{proof} 
  The lower bound when $0<\alpha<1$ follows from Lemma \ref{tva}. The upper bound follows in the same way as in the proof of Theorem  \ref{u77}. The lower bound  when $\alpha>1$ it is more complicated. For full proof see our recent paper \cite{Andersson3}.
\end{proof}

In analogy with  conjecture \ref{con3} we believe the following. 
\begin{conj} \label{con4} One can choose $B(\alpha)=\underline B(\alpha)=\overline B(\alpha)$ in Theorem \ref{avs}. \end{conj}
It should not really matter much if $|z_k| \geq 1$ or $|z_k| =1$. This has however been surprisingly difficult to prove. The technique of using Fej\'er kernels requires that $|z_k|=1$. It is possible that the method can be modified to cover the more general case, but it is not quite clear how. Nevertheless we feel  safe in believing the following conjecture.
\begin{conj} \label{con5} One has that $A(\alpha)=B(\alpha)$ where $A(\alpha)$ and $B(\alpha)$ are defined by conjectures \ref{con3} and \ref{con4}. \end{conj}
This {\em strongly} suggests that $A(\alpha)>1$ when $\alpha>1$. We can also consider (See Tur\'an \cite{Turan} page 81-83).
\begin{prob} \label{trn} Let $\alpha>1$ be a real number. Find an asymptotic formula for
\begin{gather*} 
  \inf_{z_k \in \C, |z_k| \geq 1} \max_{\nu=1,\dots,
        n^\alpha} \left| \sum_{k=1}^n z_k^\nu \right|,
 \end{gather*}
\end{prob}
From Lemma \ref{tva} and Theorem \ref{ett} we obtain similarly as in Theorem \ref{avs} that 
\begin{gather*} 
  \inf_{z_k \in \C, |z_k| = 1} \max_{\nu=1,\dots,
        n^\alpha} \left| \sum_{k=1}^n z_k^\nu \right| \sim \sqrt n, \qquad \text{for} \qquad 1<\alpha \leq 2,
 \end{gather*}
and it seems reasonable to conjecture that the same is true if $\abs{z_k} \geq 1$ instead of  $\abs{z_k} =1$.
For $\alpha>2$ the situation seemed until recently particularly unclear (see version 2 of this paper on ArXiv \cite{Andersson5}, page 19-20). However in our recent paper \cite{Andersson4} we settled an open problem of Montgomery and while we have not yet solved Problem \ref{trn} we have managed to obtain the correct order of magnitude. We proved that

 \begin{gather*} 
\sqrt{n m} \ll \inf_{z_k \in \C, |z_k| \geq 1} \max_{\nu=1,\dots,
        n^{2m}} \left| \sum_{k=1}^n z_k^\nu \right| \ll m \sqrt n.  \qquad (2 \leq m \leq n)
 \end{gather*}
 The lower bound comes from Theorem 2 in  Andersson \cite{Andersson} and the upper bound was proved by using an estimate for character sums over finite fields of Katz \cite{Katz}.

\providecommand{\bysame}{\leavevmode\hbox to3em{\hrulefill}\thinspace}
\providecommand{\MR}{\relax\ifhmode\unskip\space\fi MR }
\providecommand{\MRhref}[2]{%
  \href{http://www.ams.org/mathscinet-getitem?mr=#1}{#2}
}
\providecommand{\href}[2]{#2}
\bibliographystyle{alphaurl}

\end{document}